\documentclass[11pt]{article}
\input{epsf.sty}
\usepackage{hyperref}
\usepackage{amsmath, amssymb}
\newtheorem {thm}{Theorem}[section]
\newtheorem {prop}[thm]{Proposition}
\newtheorem {lem}[thm]{Lemma}

\def\Cox{\hfill \Box}

\def\R{{\Bbb R}}
\def\P{{\Bbb P}}

\def\0{{\bf 0}}

\def\b{\beta}

\def\d{\delta}

\def\e{\varepsilon}

\def\phi{\varphi}

\def\l{\lambda}

\def\s{\sigma}

\def\T{\T}

\def\Norm{\text{Norm}\,}

\def\Cox{\hfill \Box}

\def\R{{\Bbb R}}
\def\P{{\mathbb P}}

\def\0{{\bf 0}}

\def\RR{{\cal R}}

\def\b{\beta}

\def\d{\delta}

\def\e{\varepsilon}

\def\phi{\varphi}

\def\l{\lambda}

\def\s{\sigma}

\def\T{\T}

\def\Norm{\text{Norm}\,}

\def\Const{\text{Const}\,}

\def\one{\hbox{J}\kern-.2em\hbox{I}}

\begin{document}

\title{How non-Gibbsianness helps a 
\\ metastable Morita minimizer to provide a stable free energy}


\author{Christof K\"ulske\thanks{Research supported by the 
DFG} \footnote{TU Berlin, Institut f\"ur Mathematik, Sekr. MA 7-4, 
Strasse des 17. Juni 136, D-10623 Berlin, Germany, \texttt{kuelske@math.tu-berlin.de},}
\footnote{EURANDOM,
P.O.Box 513, 
5600 MB Eindhoven, The Netherlands
}}

\maketitle

\begin{abstract}
We analyze a simple approximation scheme based 
on the Morita-approach for  the example of the 
mean field random field Ising model 
where it is claimed to be exact in some of the physics literature. 
We show that the approximation scheme is flawed, but
it provides a set of equations whose metastable solutions surprisingly 
yield the correct solution of the model. We explain how the same equations 
appear in a different way  as rigorous consistency equations.  We 
clarify the relation between the validity of their 
solutions and the almost surely discontinuous behavior of the single-site 
conditional probabilities.

\end{abstract}

\bigskip
\bigskip


\section{Introduction} \label{sect:intro}

The {\it Morita-approach} or {\it equilibrium ensemble approach} 
to systems with quenched disorder goes back to \cite{Mor64}. 
A fair and clear review from a theoretical physicist's point of view containing a quick outline of the 
theory and various interesting 
recent applications is given by K\"uhn in \cite{Kueh04} (see also \cite{Kueh96}). 

The central idea in the Morita-approach is:   
Look at the {\it joint measure} governing the 
distribution of the quenched degrees  of freedom and the dynamical variables, 
rather than directly trying to describe 
the {\it quenched measure} for the dynamical variables 
for fixed realization of the disorder.   
Ideally one would then like to write this joint measure as a formal equilibrium 
average over the {\it joint variables} in terms of a {\it joint Hamiltonian}. This joint Hamiltonian 
would then  be the sum of the original one and another "disorder-Hamiltonian" depending 
only on the quenched degrees of freedom. The resulting model possesses full 
spatial symmetries and one might hope that it 
is amenable to techniques known from systems 
without disorder. 

Mathematically there are problems with this idea.  
In finite volume this "disorder-Hamiltonian" can in principle be chosen in such a way 
that the resulting joint distribution coincides with that of the true model. 
For lattice systems in the infinite volume this is however a serious problem. 
In fact, for many models an absolutely summable joint Hamiltonian does not exist, 
and the joint measures in the infinite volume are non-Gibbsian measures. The appearance of a non-Gibbsian 
joint distribution was first discovered in the example of the 
Grising model in \cite{vEMSS}
and studied in a general context in \cite{K1,K2}.   
See also the discussion in \cite{EKM,KM1,KuehMaz00}.
A well-understood example for this are in 
particular the joint measures of the random field Ising model \cite{BK}
in more than $3$ dimensions at low temperature and small disorder. 
They provide an illuminating example of strong non-Gibbsian 
pathologies. In fact,  their conditional probabilities are shown to be 
discontinuous functions of the conditionings, for a set of conditionings 
with (joint) measure one \cite{K2}. This means 
that the measure is not even "almost surely Gibbs" (in the sense of \cite{MRM}). This pathology 
even causes the usual Gibbs variational principle to fail \cite{KLR02}. 
Close analogies to this behavior on the lattice 
can be already found in the corresponding mean-field model. 
Here the corresponding functions describing the 
conditional expectations can be explicitly computed \cite{K3}.
For more on the analogies between non-Gibbsian measures on the lattice 
and discontinuous behavior of conditional probabilities in mean-field models 
see \cite{KuH04} and \cite{K3}.  

The motivation for using 
the Morita-approach from the point of view of theoretical physics 
is however to leave these conceptual problems aside and 
take it as a source for approximation schemes \cite{Kueh04,KM1}. 
Such schemes can be obtained 
by taking certain simplified trial disorder-Hamiltonians that 
are chosen e.g. demanding that  a finite number of moments of the 
distribution of the disorder variables coincide with that of the true distribution. 
Then one would like to solve
the resulting Morita approximant model and hope that relevant 
features of the solution are the same as that of the true model. 

It might seem hopeless to justify such 
approximations in general for non-trivial lattice models. 
It is therefore valuable to fully understand at least simple toy models that 
can be explicitly treated. This is want we want to do here. 
We will give here a complete discussion 
of the quick naive "solution" of the mean-field 
random field Ising model, based on a very simple approximant joint measure containing only 
one parameter \cite{SerPal93}.  This so-called solution is fairly old and the computations 
are trivial, but the justification of the resulting equations is subtle. So it is 
worth to reconsider it from a rigorous point of view and straighten out 
some wrong claims in the literature, answering a question of K\"uhn.  
In particular we take issue with the statement of 
K\"uhn \cite{Kueh04}, describing the work of \cite{SerPal93}.  He writes:   
"it reproduces the exact solution at the cost of introducing a single `chemical 
potential'  to  fix the average value of the random field, which creates a term in 
the modified Hamiltonian that introduces no non-locality into the model over and 
above that already contained in the definition of the Curie-Weiss limit. 
This is a remarkable result in the light of concerns raised about the 
appearance of non-Gibbsian measures within the equilibrium ensemble approach 
and the identification of the RFIM as providing a realization of a kind of 
`worse-case scenario' in the non-Gibbsian world..." 

We will indeed see that even for this simple model the situation is subtle and the validity 
of the solution is  fundamentally related to the analogue of 
"non-Gibbsianness" in the mean-field context. 
To see this, we will start in Section 2 by reviewing the quick "solution", following \cite{Kueh04}. 
This provides us with two equations for two parameters, the magnetisation and 
the Morita chemical potential. 
These equation have in fact solutions for which the magnetisation-variable 
takes the known value of the spontaneous magnetization. 
However, we note that 
this solution corresponds to a wrong (metastable) saddle-point approximation for the 
approximant measure and therefore the naive derivation given above is 
flawed.  
Moreover, we will prove that in the low temperature regime it is even strictly impossible 
to choose a chemical potential such that the magnetic fields become symmetric. In brief, 
the theory based on the Morita approximant measure with just one chemical potential fails. 

How can we understand then that the two equations derived by 
a wrong line of argument yield the correct value of the magnetization? Is this 
just accidental? 
We will see in Section 3 that the same two equations come up in a different way as 
consistency equations for the conditional probabilities of the true joint measures 
of the model without approximations.  Here however the fixed Morita-chemical potential 
is replaced by a random variable. 
It is in this context that we will finally understand that the validity of these 
equations and the almost sure discontinuity of 
the conditional expectations are consequences of each other.\bigskip 

\section{Invalidity of single-site Morita approximation 
approach for the Curie Weiss Random Field Ising Model}

\label{sect:Morita}

We consider the mean-field random field Ising model. 
It is defined in terms of the following formula for the quenched 
Gibbs expectation for fixed choice of the random fields. 
\bigskip

\noindent{\bf Quenched measure:} 
\begin{equation}\label{quenchedmeas}
\mu_{\b,\e,h_0, N}[\eta_{[1,N]}](\s_{[1, N]})
:=\frac{\exp\Bigl(\frac{\b}{2 N}\bigl(\sum_{i=1}^N \s_i\bigr)^2 + \b
 \sum_{i=1}^N (\e\eta_i+h_0)\s_i
\Bigr)}{Z_{\b,\e,h_0,N} [\eta_{[1,N]}]}
\end{equation}
Here the spins ("dynamical variables")  take values $\s_i=\pm 1 $ and 
the random fields take values $\eta_i=\pm 1$ with equal probability. 
We denote their distribution by $\P$. 
We stress that the partition function appearing in the 
denominator depends on the realization of the random 
fields $\eta_{[1,N]}$ describing the disorder. 
We allow from the beginning also an external magnetic field $h_0$, but 
we are mainly interested in the case $h_0\downarrow 0$.

What one understands by the "solution of the model" is the characterization 
of the behavior of this measure on the $\s$'s for a large set of $\eta$'s, 
having asymptotically $\P$-measure one.  This has been done in great 
detail \cite{SW85,APZ92,K97}, 
and so in this model 
there is no need for any approximation based on the Morita-approach in order to solve the model. 
Most basically, we know the phase structure in zero external field, 
for any choice of the parameters $\b,\e$.  
We recall that for large $\b$ and small $\e$ in zero 
external magnetic field $h_0$ the 
model exhibits a spontaneous magnetization whose value $m$ is a solution 
of the equation 
\begin{equation}\label{eos}
m=\frac{1}{2}\Bigr(\tanh \b(m+\e)+\tanh \b(m-\e)\Bigl)
\end{equation}

We also know finer properties of the quenched distribution 
above, like its dependence on the volume label $N$, for fixed realisation of the 
random fields. This can be asymptotically described by the in the "metastates formalism" \cite{K97}, a notion 
due to Newman and Stein \cite{NS96a,NS96b}.  
For general background on this notion in the theory of disordered systems see \cite{New97,Bo01}.

Knowing the correct solution, 
our point in this note will be however to put the Morita approximation scheme
outlined above to the test. 
Now, in the Morita-approach one looks at the joint measures on the product 
space of the spin variables $\s$ and the disorder variables $\eta$. 
They are simply composed from the quenched measures 
and the a priori uniform distribution of the random fields by the following obvious formula. 
\bigskip

\noindent{\bf True joint measure:} 
\begin{equation}
K_{\b,\e,h_0, N}(\s_{[1,N]}, \eta_{[1, N]})
=\frac{1}{2^N}\mu_{\b,\e,h_0, N}[\eta_{[1,N]}](\s_{[1, N]})
\end{equation}

The approximant measure we want to consider is obtained 
by putting a single-site disorder potential 
$\l\sum_{i=1}^N \eta_i$ with just one free parameter $\l$ that has the meaning 
of a chemical potential governing the mean value of the random fields.  
\bigskip

\noindent {\bf Morita-approximant measure:}
\begin{equation}
\hat K_{\l;\b,\e,h_0, N}(\s_{[1,N]}, \eta_{[1, N]})
=\frac{\exp\Bigl(\frac{\b}{2 N}\bigl(\sum_{i=1}^N \s_i\bigr)^2 + \b
 \sum_{i=1}^N (\e\eta_i+h_0)\s_i+\l\sum_{i=1}^N \eta_i
\Bigr)}{Z_{\l;\b,\e,h_0,N} }
\end{equation}
We stress that the partition function does {\it not} depend 
on $\eta_{[1, N]}$, in contrast to (\ref{quenchedmeas}). 
The Hamiltonian of this measure contains no non-local couplings 
of the random fields. 

Then the idea of the naive Morita approximation-approach is as follows: 
1) For any fixed $\l$, compute 
the large-$N$ limit of distribution of this model. 
2) Choose $\l=\l(\b,\e,h_0)$ such that the expectation 
of the random fields coincides with the true joint measures, i.e. it vanishes,  
$\lim_{N}\int\hat K_{\l;\b,\e,h_0, N}(d\eta_1)\eta_1=0$.
More precisely the value of $\l$ will depend on $N$, but it will 
have a well-defined limit as $N\uparrow \infty$. 
3) Then, the distribution of the Morita approximant measure  taken with 
this value of the bias of the random fields $\l$, should be close 
to the true joint measure. E.g. we should have that the distribution of the 
spin average $\frac{1}{N}\sum_{i=1}^N \s_i$ has the same infinite volume 
limit in the true joint measure and in the  Morita 
approximant measure. 

Let us write down the following precise formulation in order to have 
a well-defined starting point of discussion. 

\noindent{\bf Single-Site approximation conjecture:}   Let $\b,\e,h_0$ be fixed. 
Then the conjecture is that there is a value $\l(\b,\e,h_0)$ such that 
$\lim_{N}\hat K_{\l;\b,\e,h_0, N}(d\eta_1)\eta_1=0$ and that for 
this value we have that
\begin{equation}
\lim_{N\uparrow\infty}K_{\b,\e,h_0, N}\Bigl(\frac{1}{N}\sum_{i=1}^N \s_i \in \,\cdot\,\Bigr)
=\lim_{N\uparrow\infty}
\hat K_{\l(\b,\e,h_0);\b,\e,h_0, N}\Bigl(\frac{1}{N}\sum_{i=1}^N \s_i \in \,\cdot\,\Bigr)\end{equation}

Is this conjecture true? How does this relate to the proved a.s. discontinuity 
of the conditional expectations of the true joint measures? 
Let us review the quick derivation of the solution of the model based on this conjecture 
(we follow here  \cite{Kueh04}). 

\noindent {\bf Naive (problematic!) derivation: } 
Look at the partition function of the Morita approximant measure,   putting  $h_0=0$ from the beginning, 
and use a simple Gaussian identity (Hubbard-Stratonovitch transformation) 
to write  

\begin{equation}\label{partition-function}
Z_{\l;\b,\e,N}=2^N \int\frac{d m}{\sqrt{2\pi /(\b N)}}
\exp\Bigl(-\b N \Phi_{\l;\b,\e}(m)  \Bigr)
\end{equation}

The function appearing in the exponent is $N$-independent and 
is given below in (\ref{phifunction}) by putting $h_0=0$. 
Using the Laplace method to compute the integral 
we must have $\frac{\partial}{\partial m}\Phi_{\l;\b,\e}(m)=0$. 
This is an equation for the minimizer $m$ of the form 
\begin{equation}\label{mfeq1}\begin{split}
&m=\frac{\sum_{k=\pm 1}\sinh\bigl( 
\b(m+\e k)
\bigr) e^{\l k}}{\sum_{k=\pm 1}\cosh\bigl( 
\b(m+ \e k )
\bigr) e^{\l k}}\cr
\end{split}
\end{equation}
The parameter $\l$ is fixed such that the mean of the 
magnetic field sum divided by $N$ vanishes, i.e. $\lim_{N\uparrow\infty} \frac{\partial}{\partial \l}
\log Z_{\l;\b,\e,N}=0$. This requires at the minimizer $m$ 
that $\frac{\partial}{\partial \l}\Phi_{\l;\b,\e}(m)=0$. 
This requires that 
\begin{equation}\label{mfeq2}\begin{split}
&e^{-2 \l}=\frac{\cosh (\b(m+\e))}{\cosh(\b(m-\e))}
\end{split}
\end{equation}
The equation (\ref{mfeq2}) shows that $m$ and $\l$ are in one-to-one correspondence 
to each other. From (\ref{mfeq1}) and (\ref{mfeq2}) follows 
the well known (and correct) mean field equation 
(\ref{eos}). So it seems that 
the Morita approximation approach becomes exact in this case and we are done. 

K\"uhn writes appropriately: {\it  This result \cite{SerPal93}  - simple and reassuring as it is - 
must be regarded as remarkable in the light of concerns raised about 
the appearance of non-Gibbsian measures within the equilibrium ensemble 
approach \cite{vEMSS,K1,EKM} and the identification of the RFIM as providing a realization 
of a system exhibiting almost surely non-Gibbsian joint measures \cite{K1,K3}.}

Indeed, we note that  the "derivation" is flawed because of the following fact. 

\noindent {\bf Worrisome fact why this derivation is wrong: } Suppose that $\b>1$, $\e>0$ and $\l>0$ are fixed. 
Then the minimum of the function $m\mapsto\Phi_{\l;\b,\e}(m)$ is attained at a unique 
positive value $m^*(\l)$ (as we will see below). 
Therefore there cannot be a pair $(m^*(\l),\l)$ satisfying (\ref{mfeq2}). 
So the solution $(m,\l)$ obtained  by (\ref{mfeq1}), (\ref{mfeq2}) corresponds to 
a wrong value for the free energy. 

\noindent {\bf The remaining Morita mystery: } 
Why does the wrong minimizer  give the correct equation for the magnetization? 

It is the purpose of this note to clarify the situation. 
We will be even more general and more careful here 
and allow for a possibly non-zero external magnetic field $h_0$. 
This we do in order to investigate whether taking the limit $h_0\downarrow 0$ 
only in the end will help us to solve the problem of this approach.

We can readily solve the model for any choice of the parameters 
of inverse temperature $\b$, strength of random fields $\e$, external field $h_0$
and Morita chemical potential $\l$. 
As usual in mean field models there is convergence to (linear combinations)
of product measures over the sites $i$. Indeed, any limit measure must 
be a mixture of product measures. This is clear by de Finetti's theorem 
since the limit of exchangeable measures inherits the property of exchangeability. 

Now, solving our simple model is almost trivial when we note that by summing over the 
$\eta$ first we obtain a resulting effective Curie-Weiss Ising model with a new 
effective homogenous magnetic field acting on the $\s$'s. The computations 
are simple and will be given below for the sake of completeness. Before we do so 
let us however state the most important consequence of this in the present context. 

\begin{thm}  {\bf (Impossibility of single-site approximation of true joint measure)}\label{thm:invalidity}
Assume that $\b>1$ and $\e> 0$ are fixed. 
Then 
\begin{equation}\begin{split}
&\Bigl\{h_0\in \R, \exists \l\in \R : 
\lim_{N\uparrow\infty}
\hat K_{\l;\b,\e,h_0, N}(d \eta)=\P(d\eta)
\Bigr\}\cr
&=\R\backslash
\Bigl[-a(\b,\e),  +a(\b,\e)  \Bigr]\end{split}
\end{equation}
where $a(\b,\e)$ is strictly bigger than zero.  
Here the symbol $\lim$ denotes a weak limit. 
\end{thm}
In words the theorem  states that 
the set of external homogenous magnetic fields $h_0$ 
for which there exists a "compensating" Morita-field $\l$ that 
reproduces the neutral i.i.d. distribution for the random fields is 
bounded away from zero for any $\b>1$.  This means that 
the approximation scheme must necessarily fail in the relevant 
low temperature regime: First of all, in zero external field 
$h_0$ it is impossible to produce asymptotically symmetric i.i.d. 
random fields by an appropriate choice of $\l$. This result however, might 
not be too surprising. 
But the theorem says more: Even choosing $h_0$ strictly positive and letting it tend 
to zero afterwards won't help us. 

Having said this, it is interesting to investigate the set of parameters 
for which the distribution of random fields becomes 
neutral  i.i.d. in more detail. Let us make the following definition. 

\noindent{\bf Neutral Set: } Fix the inverse temperature $\b>0$ and $\e>0$. 
We call the parameter set 
\begin{equation}\begin{split}
&{\cal R}(\b,\e):=\Bigl\{(h_0,\l)\in \R\times \R, 
\lim_{N\uparrow\infty}
\hat K_{\l;\b,\e,h_0, N}(d \eta)=\P(d\eta)
\Bigr\}\cr
\end{split}
\end{equation}
the neutral set. Obviously ${\cal R}(\b,\e)=-{\cal R}(\b,\e)$ by the symmetry of the model. 

Then we have the following theorem.  

\begin{thm}  {\bf (Structure of neutral set)}\label{thm:structurenl}
Assume that $\b>0$ and $\e> 0$ are fixed. Then the set $\RR(\b,\e)$ is the union 
of two semi-infinite curves, related to each other by reflection
at the origin.  They are connected if and only  if $\b\leq 1$. \\

More precisely these curves are of the following form.

There is a continuous increasing function $h_0\mapsto l_{\b,\e}(h_0)$ 
that is defined an open interval of the form $(a(\b,\e),\infty)$ and takes positive values. 
The left endpoint of the interval
satisfies $a(\b,\e)\begin{cases}>0 \text{ for }  \b>1\cr =0 \text{ for }\b\leq1\end{cases}$.

Define  
\begin{equation}\label{nlgraph}\begin{split}
&\RR^+(\b,\e):=\Bigl\{
(h_0,-l_{\b,\e}(h_0))\Bigl | h_0\in (a(\b,\e),\infty)
\Bigr\}
\end{split}
\end{equation}

Then,  the neutral set has the form  
\begin{equation}\label{nldecomp}\RR(\b,\e)=\begin{cases} \RR^+(\b,\e)\cup \Bigl( -\RR^+(\b,\e) \Bigr)& \text{ for } 
\b>1 \cr \RR^+(\b,\e)\cup \Bigl( -\RR^+(\b,\e) \Bigr)\cup \Bigl\{(0,0)\Bigr\}& \text{ for }Ê\b\leq 1\end{cases}
\end{equation}

Hence, for $\b>1$ the neutral set is disconnected. 
For $\b\leq 1$ we have moreover $l^*_{\b,\e}(0+)=0$, and hence 
the neutral set is connected. 

\end{thm}

This result on the neutral set is a consequence of the solution 
of the Morita approximant for any choice of the parameters. 
We will now describe the behavior of the Morita approximant 
for general choice of the parameters. Then we will derive 
as a conclusion the explicit condition for the neutral set. 

We need some definitions.  
Define the effective magnetic field-like parameter 
\begin{equation}\begin{split}
\hat h =h_0 +\bar h_{\b,\e}(\l)
\end{split}
\end{equation}
with the function 
\begin{equation} \label{eq:10}
\bar h_{\b,\e}(\l):=\frac{1}{2\b}\log\frac{\cosh(\l+\b\e)}{\cosh(\l-\b\e)}
\end{equation}

Define the joint single-site measures depending on the parameter set, and 
on an additional (magnetization-like) parameter $m\in \R$. 
\begin{equation}\label{limitform}\begin{split}
\pi_{\l;\b,\e,h_0}[m](\s_{i},\eta_{i})&:=\frac{\exp\Bigl( \b\bigl( 
(m+\e \eta_i +h_0) \s_i
+\l\eta_i
\bigr)
\Bigr)}{2 \sum_{k=\pm 1}\cosh\bigl( 
\b(m+\e k +h_0)
\bigr) e^{\l k}
}\cr
&=
\frac{\exp\bigl(\b (m +\hat h )\s_i\bigr)}{2 \cosh 
\b (m +\hat h )}\frac{\exp\bigl((\b\e\s_i+\l) \eta_i \bigr)}{2 \cosh(\b\e\s_i+\l)}\cr
\end{split}
\end{equation}

Here we found it convenient to express the joint distribution on the r.h.s. 
appearing under the $i$-product in the form $\text{Prob}(\s_i,\eta_i)=\text{Prob}(\s_i)\text{Prob}(\eta_i|\s_i)$.
In this way the marginal on the $\s$'s can be readily read off. 
We see that the role of the parameter $\hat h$ is to 
provide an "effective magnetic field" acting on the spins.

\begin{thm}  {\bf (Solution of Morita approximant)}\label{thm:Morita approximant}
Assume that the parameters \linebreak $\b,\e\in(0,\infty)$ and $\l,h_0\in (-\infty,\infty)$ are fixed. \\

\noindent{\bf (i): }Assume at first that $\hat h\neq 0$. 
Then we have the weak convergence 
\begin{equation}\label{limitform1}\begin{split}
&\lim_{N\uparrow\infty} \hat K_{\l;\b,\e,h_0, N}(\s_{[1,N_0]},\eta_{[1,N_0]})
=\prod_{i=1}^{N_0} \pi_{\l;\b,\e,h_0}\bigl[m^{\text{CW}}(\b,\hat h)\bigr](\s_{i},\eta_{i})
\end{split}
\end{equation}
Here we have denoted by $m^{\text{CW}}(\b,h)$ 
the solution of $m=\tanh(\b(m+h))$ that has the sign of $h$, for $h\neq 0$. \\

\noindent{\bf (ii): }For $\hat h=0$ we have the weak convergence to the symmetric linear combination 
of product measures 
\begin{equation}\begin{split}
&\lim_{N\uparrow\infty} \hat K_{\l;\b,\e,h_0, N}(\s_{[1,N_0]},\eta_{[1,N_0]})\cr
&=\frac{1}{2}\prod_{i=1}^{N_0} \pi_{\l;\b,\e,h_0}\bigl[m^{\text{CW}}(\b,0+)\bigr](\s_{i},\eta_{i})
+\frac{1}{2}\prod_{i=1}^{N_0} \pi_{\l;\b,\e,h_0}\bigl[m^{\text{CW}}(\b,0-)\bigr](\s_{i},\eta_{i})
\end{split}
\end{equation}

\end{thm}

\noindent{\bf Remark: } 
Of course  
$m^{\text{CW}}(\b,h)$ is the magnetization of an ordinary Curie Weiss Ising model in an external field.

\bigskip 
\noindent {\bf Proof:}
Let us write the 
Morita approximant joint measure as a marginal  on the $\s$'s times 
the conditional measure of the random fields given the $\s$'s, that is 
\begin{equation}\begin{split}
&\hat K_{\l;\b,\e,h_0, N}(\s_{[1,N]}, \eta_{[1, N]})\cr&
=\frac{\exp\Bigl(\frac{\b}{2 N}\bigl(\sum_{i=1}^N \s_i\bigr)^2 + \b h_0
 \sum_{i=1}^N \s_i + \sum_{i=1}^N \log \cosh(\b\e\s_i+\l)
 \Bigr)}{\Norm. }
\prod_{i=1}^N \frac{e^{ (\b\e\s_i+\l) \eta_i } }{2 \cosh(\b\e\s_i+\l)}
\end{split}
\end{equation}
This shows us that the marginal distribution on the $\s$'s is given 
by an ordinary ordered mean field Ising model of the form 
$\varpropto \exp\Bigl(\frac{\b}{2 N}\bigl(\sum_{i=1}^N \s_i\bigr)^2 + \b \hat h 
 \sum_{i=1}^N \s_i  \Bigr)/\Norm. $
with the effective field $\hat h$. From here the limit statements are 
obvious by the known convergence results of the Curie Weiss Ising model 
to the corresponding (linear combination of) product measures. $\Cox$

\bigskip
Now, from the explicit solution we may derive 
explicit information on the neutral set. 
In order to do so, note at first the elementary properties 

\begin{eqnarray} \label{eq:emp-meas}
\bar h_{\b,\e}(\l)&= & \left\{ 
\begin{array}{ll}
\downarrow -\e, & \mbox{for } \l\downarrow -\infty \\
0, & \mbox{for } \l=0\\
\uparrow\e
& \mbox{for} \l\uparrow\infty \, ,
\end{array} \right.
\nonumber
\end{eqnarray}
and it is a monotonically increasing in $\l$ and odd. 
It maps $\R$ to the interval $(-\e,\e)$.

\begin{thm}  {\bf (Explicit description of neutral set)}\label{thm:explicitnl}
Assume that $\b>0$ and $\e> 0$ are fixed. The decomposition (\ref{nldecomp}) 
holds with 
\begin{equation}\begin{split}
&\RR^+(\b,\e)=\Bigl\{
(h_0,-l)\Bigl| 0<\bar h_{\b,\e}(l)< h_0, \cr
&\qquad\qquad  m^{\text{CW}}\Bigl(\b,h_0- \bar h_{\b,\e}(l)\Bigr)=\frac{\sinh( 2 l)   }{\sinh(2 \b\e) }
\Bigr\}
\end{split}
\end{equation}
This set can be written as a graph in the form 
(\ref{nlgraph}) with a continuous increasing $l_{\b,\e}(h_0)$  that 
maps the interval $(a(\b,\e),\infty)$ onto the interval $(\bar h_{\b,\e}(a(\b,\e)), \b \e)$ 
where $a(\b,\e)$ is uniquely given by 
$\bar h_{\b,\e}\bigl(a(\b,\e)\bigr)= \frac{1}{2}\sinh^{-1}\Bigl(
\sinh(2 \b \e)\,m^{\text{CW}}\bigl(\b,0+\bigr) \Bigr)$. 
\end{thm}

\noindent{\bf Remark: }Note that the above expression for $a(\b,\e)$ implies 
that $a(\b,\e)=0$ if and only if the spontaneous magnetization $m^{\text{CW}}\bigl(\b,0+\bigr)$
vanishes, i.e. $\b\leq 1$.

\bigskip 
\noindent {\bf Proof:} Suppose that $\b >1$. Then, in order 
to have convergence to a symmetric product measure 
on the random fields we must have that the parameters $\l;\b,\e,h_0$ are such that 
$\hat h \neq 0$. (Indeed, for 
$\hat h=0$ the distribution on the spins 
converges to a symmetric mixture of two different product measures. 
But from this it is obvious that also the random field distribution 
will be a mixture between two different product measures.) 
This shows that $(0,0)\not\in \RR(\b,\e)$
in that case. 

Suppose however $\b\leq 1$. Then $(h_0,\l)=(0,0)$ 
implies $\hat h=0$ which implies that the distribution 
of the $\s$'s is a symmetric product measure. But this implies 
that the distribution of the random fields will be a symmetric product measure 
so that $(0,0)\not\in \RR(\b,\e)$
in that case. \\

So, we are left with the case $\hat h\neq 0$.
We can treat the cases $\b>1 $ and $\b\leq 1$ 
on a unified basis. 
Now, conditional on the value of $\s$ the $\eta_i$ have 
an expectation value of $\tanh(\l+\b \e \s_i)$. 
We use the simple identity 
\begin{equation}\begin{split}
&\tanh(\l+\b \e \s_i)=\frac{B(1-L^2)\s_i  + L(1- B^2)}{1 - B^2 L^2}\text{ where}Ê\cr
&L=\tanh\l, \quad B=\tanh \b \e
\end{split}
\end{equation}
for $\s_i=\pm 1$. 
So the distribution on the random fields $\eta_i$ converges 
weakly to a product measure with individual expectation value 
\begin{equation}\begin{split}
&\lim_{N\uparrow\infty}
\int\hat K_{\l;\b,\e,h_0, N}(d \eta_1)\eta_1=\frac{B(1-L^2) m^{\text{CW}}(\b,\hat h) + L(1- B^2)}{1 - B^2 L^2}
\end{split}
\end{equation}

Put $l=-\l$ and use $\tanh(l)/(1-\tanh^2(l))= \sinh (2l)$. 
So, in order to have the desired convergence to the symmetric 
product measure we must have 
\begin{equation}\label{consistency}m^{\text{CW}}(\b,h_0-\bar h_{\b,\e}(l)) =\frac{\sinh(2 l)}{\sinh(2\b\e)}
\end{equation} 
This equation can only hold if $h_0-\bar h_{\b,\e}(l)$ and $l$ have the same 
sign.  By symmetry we can assume that $l >0$. But this implies that $h_0>0$ (
since $\bar h_{\b,\e}(l)>0$.) 

So, it suffices to look for all pairs $(l,h_0)$  with $l>0$ that satisfy the consistency equations
(\ref{consistency}).
The small trick we are using now is to fix the $l$ and ask for $h_0$ rather 
than doing it the opposite way. 
Fixing $l$ we see that the l.h.s. runs monotonically through the open interval 
$(m^{\text{CW}}(\b,0+),1)$ when $h_0$ runs in the "allowed range" 
$(\bar h_{\b,\e}(l),\infty)$. 

So, the set of  $l>0$ such that there exists a solution $h_0$ is 
determined by the condition $(m^{\text{CW}}(\b,0+),1)\ni  \frac{\sinh(2 l)}{\sinh(2\b\e)}$. 
Equivalently, this is the open interval $l\in (a(\b,\e),\b \e)$. Moreover 
the map to $h_0$ is continuous and monotone by known properties of the function 
$m^{\text{CW}}(\b,h)$. So it can be inverted and this yields the claim. 
$\Cox$

\bigskip
\bigskip 

So what has happened in the naive (but wrong) derivation of the mean-field equations
(\ref{mfeq1}) and  (\ref{mfeq2})? 
In order to see this let us write down a representation of the 
finite-$N$ approximant measures. 
As a result of a Gaussian transformation on the level 
of measures we get the following formula. 
\begin{prop}
In finite volume $N$ we have the identity 
\begin{equation}
\begin{split}
&\hat K_{\l;\b,\e,h_0, N}(\s_{[1,N]},\eta_{[1, N]})
=\frac{\int dm
\exp\bigl(   -\b N \hat \Phi_{\b,\e,h_0}(m)\bigr)}{\int d\tilde m
\exp\bigl(   -\b N \hat \Phi_{\l;\b,\e,h_0}(\tilde m) \bigr)}\prod_{i=1}^N \pi_{\l;\b,\e,h_0}[m](\s_{i},\eta_{i})
\end{split}
\end{equation}
Here 
\begin{equation}\label{phifunction}\begin{split}
 \hat \Phi_{\l;\b,\e,h_0}(m) 
& =\frac{m^2}{2}-\frac{1}{\b}
 \log \sum_{k=\pm 1}\cosh\bigl( 
\b(m+\e k +h_0)
\bigr) e^{\l k}\cr
&=\frac{m^2}{2}-\frac{1}{\b}\log\cosh \bigl( 
\b(m+\hat h)
\bigr)+\Const(\b,\e)
 \end{split}
\end{equation}
where $\Const(\b,\e)$ does not depend on $m$. 
\end{prop}

\noindent{\bf Remark: }
The second equality for $\hat\Phi$ can be seen e.g. by
 reexpressing the first $\cosh$ as a sum over a spin $s=\pm 1$ and 
exchanging the $s$ and $k$-sums.

\noindent{\bf Proof: }
We use a Gaussian transition kernel from the $\s$-variables 
to an auxiliary real valued variable $m$ given by $T(d m|\s_{[1,N]})=
\exp\Bigl(    -\frac{\b N}{2}   \Bigl(m-\frac{\sum_{i=1}^n\s_{i}}{N} \Bigr)^2\Bigr)dm /\Norm.$. 
We define a "big joint measure" on the spins, the random fields and also the 
auxiliary magnetization-like continuous variable by the formula 
\begin{equation}
\hat M_{\l;\b,\e,h_0, N}(d m,\s_{[1,N]}, \eta_{[1, N]})
:=\hat K_{\l;\b,\e,h_0, N}(\s_{[1,N]}, \eta_{[1, N]})T(d m|\s_{[1,N]})
\end{equation}
We see that $m$ concentrates very nicely around the value of 
the empirical average of the true spins in this measure. 
Then the non-normalized density of this "big joint measure" is given by
$\exp\Bigl(    -\frac{\b N}{2} m^2
+\b\sum_{i}\bigl( 
(m+\e \eta_i +h_0) \s_i
+\l\eta_i
\bigr)
\Bigr)$. Use this to express the "big joint measure" in the form of a marginal on the $m$ times a conditional 
measure on the $(\s,\eta)$ given the $m$. From here it is simple to get the desired formula. 
$\Cox$
\bigskip

So, conditional on a value of $m$, 
the pairs $(\s_i,\eta_i)$ are independent. 
We have then for their conditional mean values
\begin{equation}\label{eq:9}\begin{split}
&
\sum_{\s_1=\pm}\pi_{\l;\b,\e,h_0}[m](\s_1)\s_1
=\frac{\sum_{k=\pm 1}\sinh\bigl( 
\b(m+\e k +h_0)
\bigr) e^{\l k}}{\sum_{k=\pm 1}\cosh\bigl( 
\b(m+ \e k +h_0)
\bigr) e^{\l k}}\cr
&
\sum_{\eta_1=\pm}\pi_{\l;\b,\e,h_0}[m](\eta_1)\eta_1
=\frac{\sum_{k=\pm 1}k\cosh\bigl( 
\b(m+\e k +h_0)
\bigr) e^{\l k}}{\sum_{k=\pm 1}\cosh\bigl( 
\b(m+\e k +h_0)
\bigr) e^{\l k}}\cr
\end{split}
\end{equation}
We remark that, with this notation, we have that (the version for general $h_0$ of) 
the saddle point equation 
(\ref{mfeq1}) is equivalent to the consistency equation for the magnetization written as 
\begin{equation}\label{mfeq1a}\begin{split}
&
m=\sum_{\s_1=\pm}\pi_{\l;\b,\e,h_0}[m](\s_1)\s_1
\end{split}
\end{equation}
The  (version for general $h_0$ of) the neutrality equation (\ref{mfeq2}) is written as 
is equivalent to 
\begin{equation}\label{mfeq2a}\begin{split}
&
0=\sum_{\eta_1=\pm}\pi_{\l;\b,\e,h_0}[m](\eta_1)\eta_1
\end{split}
\end{equation}

Now, the large-$N$ limit of the model is obtained by looking at the absolute minimizer 
of the function $m\mapsto  \hat \Phi(m) $.  
But note that the representation for $\hat\Phi(m)$ given in the second 
line shows that is has the double-well form of the corresponding function in 
an Ising model in the external field $\hat h$. 
It is an elementary property of this function 
that its absolute minimizer has the same sign as $\hat h$. 
But this shows that the relation (\ref{mfeq2}) can not be true for 
the absolute minimizer. Instead the solution of (\ref{mfeq1}), (\ref{mfeq2}) corresponds 
to the second local minimum which is not the absolute minimum but the metastable minimum.

\bigskip
\bigskip

\section{Validity of consistency equations and almost sure discontinuity of conditional expectations} 
\label{sect:true Morita}

So how can we understand the fact that the correct 
solution of the model is obtained by 
solving equations (\ref{mfeq1a}) and (\ref{mfeq2a}) although the solution 
corresponds to the wrong saddle point? The solution to this 
puzzle is due to the fact that the naive equations have rigorous counterparts in the following sense. 
The equations we are going to state now appear as consistency equations for 
the conditional probabilities of the {\it true} joint measures. 

\begin{prop}{\bf (Consistency equations for true joint measure)}
There is a function $\l_N(\eta_{[2,N]})$, depending on the parameters
$\b,\e,h_0$, which is invariant under permutation of $(\eta_i)_{i=2,\dots,N}$
 such that we have 
\begin{equation}\begin{split}\label{true1}
&\sum_{\s_1}
K_{\b,\e,h_0, N}(\s_1)\s_1\cr
&= \sum_{\s_{[2,N]}, \eta_{[2, N]}}K_{\b,\e,h_0, N}(\s_{[2,N]}, \eta_{[2, N]})
\,\,\Biggl(\sum_{\s_1=\pm}\pi_{\l_N(\eta_{[2,N]});\b,\e,h_0}\Bigl[\frac{1}{N}\sum_{i=2}^N \s_i  \Bigr](\s_1)\s_1\Biggr)
\end{split}
\end{equation}
\begin{equation}\begin{split}\label{true2}
&0=
 \sum_{\s_{[2,N]}, \eta_{[2, N]}}
K_{\b,\e,h_0, N}(\s_{[2,N]}, \eta_{[2, N]})\,\,\Biggl(\sum_{\eta_1=\pm 1}
\pi_{\l_N(\eta_{[2,N]});\b,\e,h_0}\Bigl[\frac{1}{N}\sum_{i=2}^N \s_i  \Bigr](\eta_1)\eta_1\Biggr)\cr
\end{split}
\end{equation}
\end{prop}

\bigskip

\noindent{\bf Proof of the proposition: }
The proof is based on the following lemma. 

\begin{lem}{\bf (Representation of conditional probability of true joint measure)}
The single-site conditional probabilities can be written in the form 
\begin{equation}\begin{split}
&K_{\b,\e,h_0, N}(\s_1,\eta_1|\s_{[2,N]}, \eta_{[2, N]})
=\pi_{\l_N(\eta_{[2,N]});\b,\e,h_0}\Bigl[\frac{1}{N}\sum_{i=2}^N \s_i  \Bigr](\s_1,\eta_1)
\end{split}
\end{equation}
where
\begin{equation}\begin{split}
\l_{N}(\eta_{[2,N]})
&=\frac{1}{2}\log\frac{Z_{\b,\e,h_0,N} [\eta_1=-,
\eta_{[2,N]}]}{Z_{\b,\e,h_0,N} [\eta_1=+,
\eta_{[2,N]}]}\cr
\end{split}
\end{equation}

\end{lem}

\noindent{\bf Proof of the lemma: }
By a simple computation we have for the single-site distribution 
\begin{equation}\begin{split}
&K_{\b,\e,h_0, N}(\s_1,\eta_1|\s_{[2,N]}, \eta_{[2, N]})\cr
&= \frac{1}{\Norm}
\exp \Bigl(\b\bigl(\frac{1}{N}\sum_{i=2}^N \s_i+\e\eta_1+h_0 \bigr)\s_1
+\frac{1}{2}\log\frac{Z_{\b,\e,h_0,N} [\eta_1=-,
\eta_{[2,N]}]}{Z_{\b,\e,h_0,N} [\eta_1=+,
\eta_{[2,N]}]}\times\eta_1
\Bigr)
\end{split}
\end{equation}
and this shows the claim. 
$\Cox$

Continuing with the proof of the proposition we use the formula for the 
conditional probabilities writing 
\begin{equation}\begin{split}
&K_{\b,\e,h_0, N}(\s_1,\eta_1)\cr
&= \sum_{\s_{[2,N]}, \eta_{[2, N]}}
K_{\b,\e,h_0, N}(\s_{[2,N]}, \eta_{[2, N]})\,\,\,\, 
\pi_{\l_N(\eta_{[2,N]});\b,\e,h_0}\Bigl[\frac{1}{N}\sum_{i=2}^N \s_i  \Bigr](\s_1,\eta_1)
\cr
\end{split}
\end{equation}
But this equation gives the equation for the magnetization (\ref{true1}) by summing over $\s_1$.  
Using the symmetry of the distribution of $\eta_1$ we get (\ref{true2}).  
$\Cox$

\bigskip
\bigskip

Let us now summarize what we know by the rigorous 
solution of the random field model about 
the limiting distribution of the pair of random quantities entering the single-site kernel $\pi$. 
In words, the distribution becomes sharp in the case of non-zero external field. 
It becomes sharp but double valued in the case of vanishing external field. 
In view of the last lemma this statement is a different way of saying that there 
is a jump of the conditional probabilities when the empirical 
random field sum of the conditioning is infinitesimally perturbed around 
its typical value $0$. Now, the rigorous statement is as follows. 
\bigskip
\bigskip

\begin{thm}{\bf (Convergence of true joint measures)}

\noindent{\bf (i)} Suppose that $h_0>0$. Then we have the weak limit 
$$\lim_{N\uparrow\infty}K_{\b,\e,h_0, N}\Biggl (\frac{1}{N}\sum_{i=2}^N \s_i \in\,\cdot\,\,\quad ,\,
\l_N(\eta_{[2, N]})\in \,\cdot \,\Biggr) 
\rightarrow \d_{m^*(h_0)}\times\d_{\l^*(h_0)}
$$
Here $(m^*(h_0),\l^*(h_0))$ is a 
solution of the consistency equations 
(\ref{mfeq1a}) and (\ref{mfeq2a}). 
\medskip

\noindent{\bf (ii)} Suppose that $h_0=0$. Then 
$$\lim_{N\uparrow\infty}
K_{\b,\e,h_0=0, N}\Biggl (
\frac{1}{N}\sum_{i=2}^N \s_i \in\,\cdot\,\,\quad ,\,
\l_N(\eta_{[2, N]})\in \,\cdot \,
\Biggl| \sum_{i=1}^N\eta_i>0\Biggr) 
\rightarrow \d_{m^*}\times\d_{\l^*}
$$
where $(m^*,\l^*)$ is the
unique solution of the consistency equations 
(\ref{mfeq1a}) and (\ref{mfeq2a}) with $m^*>0$ (and, as a consequence $\l^*<0$). 

As a consequence we have 
$$
\lim_{N\uparrow\infty}K_{\b,\e,h_0=0, N}\Biggl (
\frac{1}{N}\sum_{i=2}^N \s_i \in\,\cdot\,\,\quad ,\,
\l_N(\eta_{[2, N]})\in \,\cdot \,\Biggr) 
\rightarrow \frac{1}{2}\d_{m^*}\times\d_{\l^*}+\frac{1}{2}\d_{-m^*}\times\d_{-\l^*}
$$
\end{thm}

\noindent{\bf Remark: }
We see that the system chooses the particular value of $\l_N(\eta_{[2, N]})$ (that has the opposite 
sign of the magnetisation) itself! 
\bigskip 

\noindent{\bf Proof: } We only sketch the proof. 
We rewrite the quotient of partition functions appearing in the definition of $\l_N(\eta_{[2,N]})$ in the form 
\begin{equation}\begin{split}\label{lam}
\l_{N}(\eta_{[2,N]})=\frac{1}{2}\log \mu_{\b,\e,h_0,N} [\eta_1=+,
\eta_{[2,N]}]\Bigl(\exp\bigl( - 2 \b \e \s_1  \bigr) \Bigr)
\end{split}
\end{equation}
From here Theorem 3.3 follows from the work done for the quenched 
model in \cite{K97,K3}. Let us focus here only on the interesting case of vanishing 
external magnetic field $h_0=0$. 
In this case it was shown that,  under the condition of positive sum of the random fields 
the empirical average of the spins concentrates sharply around 
the positive magnetisation $m^*$ (positive solution of (\ref{eos}))
w.r.t. to the quenched Gibbs probability. 
(This is true for "typical values" of the random field sum, 
that is for $N^{\frac{1}{2}-\d}\leq 
\sum_{i=1}^N\eta_i \leq N^{\frac{1}{2}+\d}$, and these values get all mass w.r.t. $\P$ in the large-$N$ limit).  
At the same time the 
quenched Gibbs probability  $\mu_{\b,\e,h_0,N} [\eta_1=+,
\eta_{[2,N]}](\s_1=+)$ aquires a sharp value that is related in a simple way to $m^*$. 
From (\ref{lam}) this gives the value of $\l^*$.
$\Cox$

\bigskip
\bigskip

Not assuming the knowledge of the solution of the quenched model 
we can reverse the argument in the following way in order to solve the model. 
Look at the consistency equations for the true joint measure (\ref{true1}),(\ref{true2}). 
Take $h_0>0$. Then it is very plausible without much a priori knowledge that the 
distribution of the 
pair $\Bigl( \frac{1}{N}\sum_{i=2}^N \s_i,\l_{N}(\eta_{[2,N]})
\Bigr)$ under the true joint measure should converge to a Dirac measure $\d_{m,\l}$. 
(This is in particular clear, if we assume that $\l_N$ has the form (\ref{lam}) and assume 
that the quenched magnetization becomes sharp for typical realization of the random fields in a 
positive homogeneous external field.)
But this means that the outer integrals in the 
rigorous consistency equations become sharp. 
So, the limiting value $(m,\l)$ must necessarily satisfy the 
naive consistency equations 
(\ref{mfeq1a}) and (\ref{mfeq2a}). 
These equations can then be solved and afterwards we let the external magnetic field 
$h_0$ tend to zero from above to discover 
the known solutions for the model. 

Let us finally see that, in the case of $h_0=0$ the 
validity of the naive equations implies that  
{\it there must be} discontinuous behavior of the conditional expectations as a function 
of the average of the random fields appearing in the conditioning. 
Indeed, suppose that $\l_N(\eta_{[2,N]})$ were a continuous function 
of $\frac{1}{N}\sum_{i=2}^N \eta_i$. Then, by the law 
of large numbers it would have to be constant in the large-$N$ limit. 
But by reasons of symmetry this constant would have 
to be zero in the case of $h_0=0$.  
But this is in contradiction to the non-trivial solution of the 
naive equations (\ref{mfeq1}) and (\ref{mfeq2}). 
To summarize the last line of 
argument in catchy terms: Non-Gibbsianness is necessary 
to help the metastable solution to provide the right answer.  

\bigskip 
We remark that the purpose of this note is not 
to attack the Morita approach in general as a valuable heuristic method 
in theoretical physics to predict the behavior of disordered systems 
when a rigorous analysis is not available or not yet available.   

As pointed out to us by Reimer K\"uhn, one could also argue that the second of
the naive equations (\ref{mfeq2}), which
demands that $m$ and $\lambda$ at the physical fixed point
must have opposite sign,  renders the region of integration for the partition function
(\ref{partition-function})
which includes the other fixed point as {\it  unphysical} and so
to be excluded from the domain of integration. 
This line of reasoning would render the "naive" argument
correct and this is not the first occasion in physics where such things
happen. 
While there seems no direct rigorous justification for this procedure  
we have shown that one is able to  understand the validity of the naive equations 
by viewing the parameter $\l$ properly as a stochastic quantity. 
This might give hope that results 
obtained by approximation schemes based on the Morita approach  
provide correct answers also in more complicated situations where 
a rigorous analysis is lacking. A better understanding of this 
would pose a fascinating challenge.

\bigskip
\bigskip

{\bf Acknowledgments:} 

The author thanks Aernout van Enter and Reimer K\"uhn for discussions at the conference 
"Gibbs vs. non-Gibbs"  (EURANDOM, december 2003) 
which led to this paper and valuable comments on an earlier draft.

\end{document}